\documentclass[11pt]{article}
\author{Bing-Long Chen and   Le Yin\\[8pt]}
\title{\textbf{Isometric embedding of negatively curved complete surfaces in Lorentz-Minkowski space}}
\date{October 16, 2011}
\usepackage{latexsym}
\usepackage{amsmath}
\usepackage{amssymb,amscd}
\usepackage[mathscr]{eucal}
\newtheorem{thm}{Theorem}[section]
\newtheorem{cor}[thm]{Corollary}
\newtheorem{lem}[thm]{Lemma}
\newtheorem{prop}[thm]{Proposition}

\newtheorem{rem}{Remark}[section]

\numberwithin{equation}{section}
\newenvironment{pf}{{\noindent \it  Proof.}}{{\hfill$\Box$}\\}
\begin{document}
\maketitle
 \let\thefootnote\relax\footnotetext{The first author was supported partially by grants NSFC11025107,
NSFC10831008, FRFCU20\\10-34000-3162643, HLTP34000-5221001, and the
second author(corresponding author)  by NSFC 11101289. }
\begin{abstract}
 Hilbert-Efimov theorem  states that any complete surface with curvature bounded above by a negative constant can not be
isometrically imbedded in $\mathbb{R}^3.$ We demonstrate that any
simply-connected smooth complete surface with curvature bounded
above by a negative constant admits a smooth isometric embedding into the
Lorentz-Minkowski space $\mathbb{R}^{2,1}$.

\end{abstract}


  \section{Introduction}
 \qquad Many closed convex  surfaces can be easily found
in $\mathbb{R}^3.$ In 1906, Weyl\cite{W} posed the problem whether
an abstract compact smooth simply-connected two dimensional
Riemannian manifold with positive curvature could be realized in
$\mathbb{R}^3.$ The problem, named after Weyl,  was investigated  by
Weyl-Lewy-Alexandrov, and finally resolved (in smooth category) by
Nirenberg \cite{Ni} and Pogorelov \cite{Po} independently.

The generalization to nonnegative curvature case was done by Guan-Li
\cite{GL} and Hong-Zuilly \cite{HZ}, and only $C^{1,1}$ imbedding
was obtained. For local isometric embeddings, there were important
breakthroughs  of Lin \cite{Lin1} \cite{Lin2}, Han-Hong-Lin
\cite{HHL}, Han \cite{Haq} \cite{Haq1}. The reader is refereed to
the survey articles  Hong \cite{Ho}, Yau \cite{Yau} and the book
Han-Hong \cite{HH}.

 The story is completely different for surfaces
with negative curvature, the famous Hilbert-Efimov
theorem(\cite{Hil},\cite{EF}) asserts that any  complete surface
with curvature bounded above by a negative constant can not be
realized in $\mathbb{R}^3.$

On the other hand, the hyperbolic plane $\mathbb{H}^2$ admits a
canonical smooth isometric embedding in the 3-dimensional
Lorentz-Minkowski space $\mathbb{R}^{2,1}$ as a unit imaginary
sphere $x_3^2-(x_1^2+x_2^2)=1.$   Here $\mathbb{R}^{2,1}$ is
$\mathbb{R}^{3}$ equipped with metric $ds^2=dx_1^2+dx_2^2-dx_3^2.$
Instead of the Euclidean space $\mathbb{R}^3$, it is proved that the
Lorentz-Minkowski space $\mathbb{R}^{2,1}$ is the appropriate
ambient space for the isometric imbedding of strongly negatively
curved surfaces.

 The problem of isometric embedding of surfaces with boundary into
$\mathbb{R}^{2,1}$ actually has been studied by B. Guan \cite{G}.
The author proved that the smooth metric of negative curvature on
2-disc ${\mathcal{D}}$ with boundary $\partial \mathcal{D}$ of
positive geodesic curvature admits a smooth isometric embedding $X:
\mathcal{D}\rightarrow \mathbb{R}^{2,1}$ with planar boundary
$X(\partial \mathcal{D})\subset \mathbb{R}^{2}.$  The purpose of the
paper is to find global isometric embeddings for complete surfaces.
The main result is the following:
\begin{thm}\label{t1} Let $(M,g)$ be a smooth two-dimensional simply-connected complete  Riemannian
 manifold with curvature $K$ satisfying
 \begin{equation} K\leq -C_1
 \end{equation}
 for some positive constant $C_1>0.$
 Then there exists a smooth isometric embedding $X: M\rightarrow
 \mathbb{R}^{2,1},$ and the spacelike submanifold $X(M)$ is a
 graph over
 $\mathbb{R}^2\subset \mathbb{R}^{2,1}:(x_1,x_2,0)\rightarrow (x_1,x_2,Z(x_1,x_2))$
 satisfying
 \begin{equation}\label{pinch-20}\sqrt{x_1^2+x_2^2}\leq Z(x_1,x_2)\leq
 \sqrt{\frac{1}{C_1}+x_1^2+x_2^2}.
 \end{equation}

   \end{thm}
\begin{rem}\label{rm1}

 It is likely that the solution of the isometric embedding problem is far from
being unique if we drop the restriction (\ref{pinch-20})(or
(\ref{pinch-200}) below).  Actually, a remarkable fact (see \cite{GJS})
is that there are many distinct isometric embeddings for the
hyperbolic plane $\mathbb{H}^2$ into $\mathbb{R}^{2,1},$ some of
them even have unbounded second fundamental forms and violate
(\ref{pinch-20})(or (\ref{pinch-200})). In this sense, the isometric
embedding provided by Theorem \ref{t1} or the following Theorem
\ref{t1.2} is rather special. The construction and classification of
all other exotic embeddings deserve a further study.
   \end{rem}

The following Theorem \ref{t1.2} refines the previous one, it handles the
estimate of extrinsic geometries of the embedding. 
The estimate is not trivial, it holds for the particular embedding
constructed in the theorem. Note that by Remark \ref{rm1}, it is
possible that some exotic embeddings may violate the estimate.
 Theorem \ref{t1.3}  asserts the uniqueness
of these particular embeddings. We remark that the boundedness of the second fundamental form is not sufficient
to guarantee the uniqueness of the isometric embedding.

\begin{thm}\label{t1.2} Let $(M,g)$ be a smooth two-dimensional simply-connected complete  Riemannian
 manifold, whose Gauss curvature satisfies  \begin{equation}\label{1.3}-{C}_2\leq K\leq
 -C_1,\end{equation}
 \begin{equation}\label{1.4}\sup_{d(x,y)\leq 1}\frac{|K(x)-K(y)|}{d(x,y)^{\mu}}\leq {{C}}_{\mu},\end{equation}
 for some positive constants $C_2\geq C_1>0$, $1>\mu>0$, ${C}_{\mu}>0.$

 Then there exists a  smooth isometric embedding $X: M\rightarrow
 \mathbb{R}^{2,1}$ such that the spacelike submanifold $X(M)$ is a
 graph over
 $\mathbb{R}^2\subset \mathbb{R}^{2,1}:(x_1,x_2,0)\rightarrow (x_1,x_2,Z(x_1,x_2))$
 satisfying

 i) \begin{equation}\label{pinch1}\sqrt{\frac{1}{{C}_2}+x_1^2+x_2^2}\leq Z(x_1,x_2)\leq
 \sqrt{\frac{1}{C_1}+x_1^2+x_2^2};
 \end{equation}

  ii)  $|A|\leq C$, where $A$ is the second fundamental form of the
  submanifold $X(M),$  constant $C$ only depends on $C_1$, ${C}_2$ and ${{C}}_{\mu}.$


   \end{thm}
\begin{rem}\label{r1.2}
If the  curvature covariant derivatives up to order $l$ are assumed
to be bounded
 in Theorem \ref{t1.2}, i.e.  \begin{equation}\label{42}\sum_{p=0}^{l}|\nabla ^pK|\leq \bar{C}_l\end{equation} for
 some $l\geq 1,$ then the covariant derivatives of the second
 fundamental form of $X(M)$ up to order $l-1$ are also bounded
 \begin{equation}\label{010}\sup_{x\in X(M)}\sum_{p=0}^{l-1}|\nabla^p A|(x)\leq C
 \end{equation}
for some $C$ depending only on
$\bar{C}_l.$
\end{rem}
\begin{thm}\label{t1.3}
Under assumptions of Theorem \ref{t1.2}, let $X$ be the isometric imbedding constructed in Theorem \ref{t1.2}. Then

 i) let $\tilde{X}$ be  another  isometric embedding of $(M,g)$ into $\mathbb{R}^{2,1}$ such that
  $\tilde{X}(M)$
 is represented as an graph
 \begin{equation}\label{pinch-200}\sqrt{y_1^2+y_2^2}\leq \tilde{Z}(y_1,y_2)\leq
 \sqrt{\frac{1}{C}+y_1^2+y_2^2},
 \end{equation}
 in some Lorentz-Minkowski coordinate system $\{y_1,y_2,y_3\},$
    then there is an isometry $\iota\in Iso(\mathbb{R}^{2,1})$ such that
 $\tilde{X}=\iota \circ X;$

 ii) there is an injective homomorphism
 $\rho: Iso(M,g) \rightarrow Iso(\mathbb{H})\subset Iso(\mathbb{R}^{2,1})$ such that
 $$
 X\circ \gamma=\rho(\gamma)\circ X $$
 for any $\gamma\in Iso(M,g),$
 where $Iso(M,g),$ $Iso(\mathbb{H}),$  $Iso(\mathbb{R}^{2,1})$ are the groups of
 isometries of $M,$ the unit imaginary  sphere in $\mathbb{R}^{2,1}$ and $\mathbb{R}^{2,1}$ respectively.

 \end{thm}

We discuss the proof of the Theorems. It is reduced to solving
certain equation of Monge-Amp\`{e}re type:
\begin{equation}\label{4}
\frac{\det(\nabla^2 u+g)}{\det(g)}=-K_g(|\nabla u|^2+2u)
\end{equation}
on the whole manifold $M.$ The corresponding Dirichlet
problems may be solved on a sequence of exhausting domains $\Omega_{l}$  with some
particular boundary values. The problem  amounts to derive certain
uniform a priori estimates for these solutions $u_l.$ The bulk of
the present paper is devoted to these estimates. 
Historically, there were two sorts of methods to derive the
second order derivative estimate (for Weyl problem or Minkowski problem).
The first one was developed by Weyl, Nirenberg, Pogorelov, Cheng-Yau etc.,
 the second one  was done by Lewy and Heinz. The former is straightforward and works for higher dimensions,
  but the argument is
hardly to be localized. The latter one is  complex in nature, but
the estimate is purely local estimate, although  it is hardly to be
generalized to higher dimensions.

The structure of the paper is the following.
 In section 2, we outline a sketch of the proof and derive the
zero and first order estimates. In section 3, we derive the second
and higher order estimates, and Theorem \ref{t1} is proved in this section.
The proof of Theorems \ref{t1.2}, \ref{t1.3} will be given in section 4.  In the appendix,
we supply an alternative,  straightforward argument for the second order derivative estimate.  
\qquad

\vskip 0.3cm
\section{Zero and first order estimates }
\subsection{Sketch of proof} Suppose $X:M\rightarrow \mathbb{R}^{2,1}$ is an
isometric embedding, then $X(M)$ is a spacelike submanifold and the
Gauss-Codazzi-Weingarten equations read as follows
\begin{equation} \label{1}
\begin{split}
&\nabla_i\nabla_j X=h_{ij}\vec{n},  \\
&\nabla_{i}\vec{n}=h_{ij}g^{jk}X_k ,\\
& R_{ijkl}=-(h_{ik}h_{jl}-h_{il}h_{jk}),\\
& \nabla_{i}h_{jk}=\nabla_{j}h_{ik},
\end{split}
\end{equation}
where $\vec{n},$ $h_{ij},$ $R_{ijkl}$ are the normal vector, second
fundamental form and the curvature tensor respectively.

Let $u=-\frac{1}{2}\langle X, X \rangle$, where
$\langle\cdot,\cdot\rangle$ is the Lorentz-Minkowski metric. By
(\ref{1}), we have
\begin{equation}\label{2}\begin{split}
\nabla_iu&=-\langle X,X_i\rangle,\\
\nabla_i\nabla_ju&=-h_{ij}\langle \vec{n},X\rangle-g_{ij}.
\end{split}
\end{equation}
 Since
\begin{equation*}
\begin{split}
   \langle X,X\rangle
   =&\sum_{i=1}^{2}g^{ij}\langle X, X_i\rangle\langle X, X_j\rangle-\langle X,\vec{n}\rangle^2\\
   =&|\nabla u|^2-\langle X,\vec{n}\rangle^2,
  \end{split}
  \end{equation*}
  then
\begin{equation}\label{3}
\langle X,\vec{n}\rangle^2=|\nabla u|^2+2u.
\end{equation}
Combining (\ref{1}) (\ref{2}) and (\ref{3}), we get
\begin{equation}\label{4}
\frac{\det(\nabla^2u+g)}{\det(g)}=-K_g(|\nabla u|^2+2u).
\end{equation}
Note that equation (\ref{4}) satisfied by the function
$-\frac{1}{2}\langle X, X \rangle$ is an intrinsic one defined on
$(M,g).$

Conversely, if we can solve equation (\ref{4}), and get a bounded
positive
 solution $u$ of (\ref{4}) on $M$,  we will show in the following that this  yields an
isometric embedding  $X:(M,g)\rightarrow \mathbb{R}^{2,1}$ such that
$-\frac{1}{2}\langle X, X \rangle=u.$

 To construct this isometric
embedding, we need to introduce the polar coordinates in the open
future timelike cone $\mathcal{I}^{+}=\{(x_1,x_2,x_3)\in
\mathbb{R}^{2,1}|\sqrt{x_1^2+x_2^2}< x_3\}.$ In this polar
coordinate system, the Lorentz-Minkowski metric takes the form
\begin{equation}\label{5}
-dr^2+r^2ds_{\mathbb{H}}^2,
\end{equation}
 where
$r=\sqrt{x_3^2-x_1^2-x_2^2}$ and $ds_{\mathbb{H}}^2$ is the
hyperbolic metric($K=-1$) of the unit imaginary sphere: $r=1$.

\begin{prop} \label{p2.1} For a positive $C^2$ function $u$ on $M,$ define a new metric
$$\bar{g}=\frac{g+(d\sqrt{2u})^2}{2u}$$
 on $M.$  The Gauss curvature $K_{\bar{g}}$ of
$\bar{g}$ is given by
\begin{equation}
K_{\bar{g}}=-1+\frac{\frac{\det{(\nabla^2
u+g)}}{\det({g})}+K_g(|\nabla u|^2+2u)}{(1+\frac{|\nabla
u|^2}{2u})^2}.
\end{equation}
\end{prop}
\begin{pf} The Gauss curvature of the
metric $g_1\triangleq g+(d\sqrt{2u})^2$ can be computed by the
formula (see \cite{G}, \cite{Ho})
$$K_{g_1}=\frac{1}{1+|\nabla \sqrt{2u}|^2}(K_g+\frac{\det(\nabla^2\sqrt{2u})}{\det(g)(1+|\nabla\sqrt{2u}|^2)}).$$
From the curvature formula of conformal transformation
$\bar{g}=\frac{g_1}{2u},$ a straightforward computation shows
\begin{equation*}
\begin{split}
 \frac{K_{\bar{g}}}{2u}=&K_{g_1}+\frac{1}{2}\triangle_{g_1}\log u\\
 =&-\frac{1}{2u}+\frac{1}{2u}\frac{\frac{\det{(\nabla^2
u+g)}}{\det({g})}+K_g(|\nabla u|^2+2u)}{(1+\frac{|\nabla
u|^2}{2u})^2}.
\end{split}
\end{equation*}
\end{pf}
\begin{rem}\label{rem1} If we can solve (\ref{4}), and
 $u$ is a bounded positive smooth solution,   then the metric $\bar{g}$ in Proposition \ref{p2.1} is complete and has
constant curvature $-1.$ Hence there exists an isometry
$i:(M,\bar{g})\rightarrow\mathbb{H}=\{r=1\}$ and we can construct an
 embedding $I: (M,g)\rightarrow
\mathcal{I}^{+}\subset\mathbb{R}^{2,1}$ as
$I(y)\triangleq(i(y),\sqrt{2u(y)})$ in the polar coordinate system
(\ref{5}). It is clear that
\begin{equation}\label{-7}
\begin{split}
I^*(-dr^2+r^2ds_{\mathbb{H}}^2)=&-(d\sqrt{2u})^2+2u
i^*(ds_{\mathbb{H}}^2)\\
=&-(d\sqrt{2u})^2+2u\bar{g}\\
=&g,
\end{split}
\end{equation}
which shows that the map $I$ is the desired isometric embedding. The
regularity of the embedding $I$ follows from the regularity of $u.$
\end{rem}

 Hence the proof of Theorem \ref{t1} may be reduced to solving
equation (\ref{4}).  The result to be proved is the
 following:

\begin{thm}\label{t2.2}
Under the assumptions of Theorem \ref{t1}, equation (\ref{4}) has a
 smooth bounded  positive solution $u$ such that $0 < u\leq
\frac{1}{2C_1}.$
\end{thm}

 The following strategy will be adapted to solve equation
(\ref{4}). We first solve the equation (\ref{4}) on a sequence of
compact smooth exhausting domains $\Omega_{1}\subset\subset
\Omega_2\subset\subset\cdots.$ Let $u_{l}$ be the solution on
$\Omega_{l}.$ Fix $x_0\in M$, we will show that for any nonnegative integer
$k\geq0,$ there exists a constant $D_{k}>0$ such that
\begin{equation}\label{est}\sup_{\Omega_{l}\supset B(x_0,k+1)}
|u_{l}|_{C^{k}(\bar{B}(x_0,k))} \leq D_{k},
\end{equation}
where the norm $C^{k}(\bar{B}(x_0,k))$ can be defined on  some (and
any) fixed finite coordinate covering of $\bar{B}(x_0,k).$ Once
(\ref{est}) has been obtained, we can extract a subsequence of $u_l$
by Arzela-Ascoli theorem such that the limit is a smooth solution of
equation (\ref{4}).

Indeed, we choose simply $\Omega_{l}=B(x_0,l)$ and consider the
Dirichlet problem
\begin{equation}\label{Diri}
\left\{
\begin{split}
\frac{\det(\nabla^2u+g)}{\det(g)}&=-K_g(|\nabla u|^2+2u).\\
u\mid_{\partial B(x_0,l)}&=\frac{1}{2C_2(l)},
\end{split}
\right.
\end{equation}
where $C_2(l)=\max\limits_{x\in \bar{B}(x_0,l)}(-K_g(x)).$

Clearly,  (\ref{Diri}) has a subsolution
$u_0\equiv\frac{1}{2C_2(l)},$ i.e. \begin{equation*}
\frac{\det(\nabla^2u_{0}+g)}{\det(g)}\geq -K_g(|\nabla u_0|^2+2u_0).\\
\end{equation*}  By continuity method, this implies that (see \cite{GB}) (\ref{Diri})
admits a smooth solution $u_l$ which satisfies $u_l\geq u_0$ and
$\nabla^2 u_l+g>0$.

The main task of the subsequent sections is to derive a priori
estimates for the solutions $u_{l}$ so that (\ref{est}) holds. For
convenience, we drop the subscript $l$ from $u_l$ and $\Omega_l$ in
the process of computations.

\subsection{Zero and first order estimates}

\begin{prop}\label{0order}
The solution $u$ of the Dirichlet problem (\ref{Diri}) satisfies
\begin{equation}
\frac{1}{2C_2(l)}\leq u\leq\frac{1}{2C_1}.
\end{equation}
\end{prop}
\begin{pf} By applying maximum principle to $u,$ we have
$$\frac{1}{2C_2(l)} \leq u\leq
\max\{\sup_{\Omega}(-\frac{1}{2K_g}),
\frac{1}{2C_2(l)}\}\leq\frac{1}{2C_1}.$$
\end{pf}
\begin{prop}\label{1order}(1st-Order Estimate) The gradient of the
solution $u$ of (\ref{Diri}) satisfies
\begin{equation}|\nabla u|\leq \frac{2}{\sqrt{C_1}}.
\end{equation}
\end{prop}
\begin{pf}
We choose
$\xi=\frac{1}{2C_2(l)}+\frac{2}{\sqrt{C_1}}(l-d(x_0,\cdot))$ as a
barrier function. Clearly $\xi$ satisfies $\xi\mid_{\partial
\Omega}=\frac{1}{2C_2(l)}$ and $|\nabla
\xi|\leq\frac{2}{\sqrt{C_1}}.$ By Hessian comparison theorem, we
have
$$\triangle\xi=-\frac{2}{\sqrt{C_1}}\triangle d\leq-2.$$

On the other hand, from $\nabla^2u+g>0,$ we know that
$$\triangle u+2>0.$$
Hence $\triangle(u-\xi)>0$ on $\Omega$. The maximum principle
implies that $u\leq\xi.$ Therefore we have
\begin{equation}\label{dxi}
|\nabla u|\mid_{\partial \Omega}\leq
|\nabla\xi|\leq\frac{2}{\sqrt{C_1}}.
\end{equation}
Now we consider the quantity $|\nabla u|^2+2u.$  The maximum
$\max\limits_{\bar{\Omega}}(|\nabla u|^2+2u)$ is achieved either on
the boundary or in the interior of the domain. In the former case,
the maximum is bounded by $\frac{4}{C_1}+\frac{1}{C_2(l)}$ by
(\ref{dxi}). In the latter case, suppose the maximum is achieved at
some point $\bar{x}\in \Omega$. Since
   \begin{equation*}
    \begin{split}
    0=\nabla_i(|\nabla u|^2+2u)(\bar{x})=2(u_{ij}+g_{ij})u_j(\bar{x})
    \end{split}
   \end{equation*}
   and $u_{ij}+g_{ij}>0$, it follows that
   $|\nabla u|(\bar{x})=0,$ and therefore
\begin{equation*}
\max_{\bar{\Omega}} (|\nabla u|^2+2u)\leq \max_{\bar{\Omega}}2u \leq
\frac{1}{C_1}.
\end{equation*}
Combining both cases, we get
$$\sup_{\bar{\Omega}}|\nabla
u|\leq\max{\{
\frac{2}{\sqrt{C_1}},\sqrt{\frac{1}{C_1}-\frac{1}{C_2(l)}}\}}=\frac{2}{\sqrt{C_1}}.$$
\end{pf}

Propositions \ref{0order} and \ref{1order} state  that the function $u$
 and its gradient can be bounded from above by a constant independent of
the domain $\Omega_l.$ Before  estimating the lower bound of $u,$
 we need to construct cutoff functions around points where the
values of $u$ are not too large.

\begin{lem}\label{39}  Fix
$\tilde{x}\in M,$ suppose there exist a real number $r_0>0$ and  a
solution $u$ of (\ref{Diri}) defined on domain  $\Omega_{l}\supset
B(\tilde{x},r_0)$ satisfying
\begin{equation}\label{assump1}
u(\tilde{x})< \frac{r_0}{2\sqrt{\max\limits_{y\in
\bar{B}(\tilde{x},r_0)}(-K_g(y))}\coth (\sqrt{\max\limits_{y\in
\bar{B}(\tilde{x},r_0)}(-K_g(y))}r_0)}.
\end{equation}
Then there are a domain
${Q}_{\tilde{x}}\subset B(\tilde{x},r_0)$ containing $\tilde{x}$ and
a function $\varphi^{\tilde{x}}\in C^{2}(\bar{{Q}}_{\tilde{x}})$
such that
 \begin{equation*}
 \begin{split}  \text{i)}\ & 0\leq \varphi^{\tilde{x}}\leq
\frac{r_0}{2\sqrt{c_2}\coth(\sqrt{c_2}r_0)},\\
&\varphi^{\tilde{x}}\geq
\frac{\frac{r_0}{2\sqrt{c_2}\coth(\sqrt{c_2}r_0)}-u(\tilde{x})}{2} \
\ \
\text{on} \ \ {B(\tilde{x},\frac{\sqrt{C_1}(\frac{r_0}{2\sqrt{c_2}\coth(\sqrt{c_2}r_0)}-u(\tilde{x}))}{6})},\\
&\varphi^{\tilde{x}}\mid_{\partial Q_{\tilde{x}}}=0,
\end{split}
\end{equation*}

ii) $ |\nabla \varphi^{\tilde{x}}|\leq \frac{3}{\sqrt{C_1}},$

iii) $\nabla^{2}\varphi^{\tilde{x}}\geq -(\nabla^2 u+g),$

\noindent where $c_2=\max\limits_{y\in
\bar{B}(\tilde{x},r_0)}(-K_g(y)).$
\end{lem}
\begin{pf} Set
$$\xi=u+\frac{d^2(\tilde{x},\cdot)}{2\sqrt{c_2}r_0\coth(\sqrt{c_2}r_0)},
\ \  Q_{\tilde{x}}=\{\xi<\frac{r_0}{2\sqrt{c_2}\coth(\sqrt{c_2}r_0)}\},$$
and
$$\varphi^{\tilde{x}}=\frac{r_0}{2\sqrt{c_2}\coth(\sqrt{c_2}r_0)}-\xi.$$
Then $\varphi^{\tilde{x}}$ satisfies i) and ii). To check that
$\varphi^{\tilde{x}}$ satisfies iii), we use the Hessian comparison
theorem
$$
\nabla^{2}d^2(\tilde{x},\cdot)\leq 2\sqrt{c_2}d(\tilde{x},\cdot)
\coth(\sqrt{c_2}d(\tilde{x},\cdot))g
$$
to conclude that
$$
\nabla^2\xi\leq  \nabla^2 u+g.
$$
The proof of Lemma \ref{39} is completed.
\end{pf}
\begin{prop}\label{0'order}(Lower bound of $u$) For any
$\tilde{x}\in M,$ $r_0>0,$ assume the solution $u$ of (\ref{Diri})
is defined on a domain $\Omega\supset B(\tilde{x},r_0),$ then  we
have
\begin{equation}\begin{split}
u(\tilde{x})\geq \min\{&
\frac{r_0}{4\sqrt{\max\limits_{\bar{B}(\tilde{x},r_0)}(-K_g(x))}\coth (\sqrt{\max\limits_{
\bar{B}(\tilde{x},r_0)}(-K_g(x))}r_0)},\\
&\frac{C_1r_0^2}{36\max\limits_{\bar{B}(\tilde{x},r_0)}(-K_g(x))\coth^2 (\sqrt{\max\limits_{
\bar{B}(\tilde{x},r_0)}(-K_g(x))}r_0)},\\
&  \frac{1}{32 \max\limits_{\bar{B}(\tilde{x},r_0)}(-K_g(x))} \}.
\end{split}
\end{equation}
\end{prop}
\begin{pf}
Assume
\begin{equation}\label{contradiction}
 u(\tilde{x})<
\frac{r_0}{4\sqrt{\max\limits_{\bar{B}(\tilde{x},r_0)}(-K_g(x))}\coth (\sqrt{\max\limits_{
\bar{B}(\tilde{x},r_0)}(-K_g(x))}r_0)}.
\end{equation}
Clearly the condition (\ref{assump1}) holds for this
$r_0.$ Consider the quantity $\frac{u}{\zeta}$ around $\tilde{x}$,
where $\zeta=\varphi^{\tilde{x}}$ be the  cutoff function in Lemma
\ref{39}.  Suppose the minimum of $\frac{u}{\zeta}$ is achieved at
some point $\bar{x}\in $ $\text{supp}(\zeta)$.  At the point
$\bar{x}$, we have
\begin{equation}\label{zeta100}
\frac{\nabla u}{u}=\frac{\nabla \zeta}{\zeta},
\end{equation}
and
\begin{equation}\label{zeta101}
0\leq\nabla^2\log\frac{u}{\zeta}=\frac{\nabla^2
u}{u}-\frac{\nabla^2\zeta}{\zeta}.
\end{equation}
Diagonalize $u_{ij}=\lambda_i\delta_{ij}$ at $\bar{x}$ with an
orthonomal basis. It follows from (\ref{zeta100}) and
(\ref{zeta101}) that
\begin{equation}\label{zeta2}
\begin{split}
\sum
\frac{\nabla_{ii}\zeta}{1+\lambda_i}
&\leq\frac{\zeta}{u}(2-\frac{2}{\sqrt{(1+\lambda_1)(1+\lambda_2)}})\\
&=\frac{2\zeta}{u}(1-\frac{1}{\sqrt{(-K_g)(|\nabla
\zeta|^2\frac{u^2}{\zeta^2}+2(\frac{u}{\zeta})\zeta)}}).
\end{split}
\end{equation}
Denote $$A=2\sqrt{\max\limits_{\bar{B}(\tilde{x},r_0)}(-K_g(x))}\coth (\sqrt{\max\limits_{
\bar{B}(\tilde{x},r_0)}(-K_g(x))}r_0).$$
Combining (\ref{zeta2}) and Lemma \ref{39}, we have
\begin{equation}\label{zeta3}
\begin{split}
-2\leq\frac{2\zeta}{u}(1-\frac{1}{\sqrt{(-K_g)(\frac{9}{C_1}\frac{u^2}{\zeta^2}
+\frac{2r_0}{A}\frac{u}{\zeta})}}).
\end{split}
\end{equation}
If $\sqrt{(-K_g)(\frac{9}{C_1}\frac{u^2}{\zeta^2}
+\frac{2r_0}{A}\frac{u}{\zeta})}\leq\frac{1}{2}$, by (\ref{zeta3})
we get
\begin{equation}\label{zeta4} \frac{u}{\zeta}\geq 1.
\end{equation}
On the other hand, if $\sqrt{(-K_g)(\frac{9}{C_1}\frac{u^2}{\zeta^2}
+\frac{2r_0}{A}\frac{u}{\zeta})}>\frac{1}{2}$, direct computation
shows
\begin{equation}\label{zeta5}
\frac{u}{\zeta}\geq \min\{\frac{A}{16r_0 \max\limits_{x\in
\bar{B}(\tilde{x},r_0)}(-K_g(x))}, \frac{2C_1r_0}{9A}\}.
\end{equation}Combining (\ref{zeta4}) and (\ref{zeta5}),  we have
\begin{equation*}
\frac{u}{\zeta}\geq \min\{1,\frac{A}{16r_0 \max\limits_{x\in
\bar{B}(\tilde{x},r_0)}(-K_g(x))}, \frac{2C_1r_0}{9A}\}.
\end{equation*}
In particular, this implies
\begin{equation*}\begin{split} u(\tilde{x})\geq
&\min\{\frac{r_0}{2A},\frac{1}{32\max\limits_{x\in
\bar{B}(\tilde{x},r_0)}(-K_g(x))}, \frac{C_1r_0^2}{9A^2}\}.
\end{split}
\end{equation*}
The proof of Proposition \ref{0'order} is completed.
\end{pf}
\begin{cor}\label{0''order} For any $r_0>0,$ there is a constant $C$
depending only on $r_0$ and $C_1$ such that
\begin{equation}\begin{split}
u(\tilde{x})\geq
\frac{C^{-1}}{\max\limits_{\bar{B}(\tilde{x},r_0)}(-K_g(x))},
\end{split}
\end{equation}
for any solution $u$ to (\ref{Diri}) defined on $\Omega\supset
B(\tilde{x},r_0).$
\end{cor}

\vskip 0.3cm

\section{Second and higher order estimates}

In this section, we will give a purely local second order derivative estimate.
This estimate could be done by Heinz-Lewy ``characteristic'' 
theory for Monge-Amp\`{e}re equations in  dimension two.
The reader is referred to the lecture notes \cite{Sch} by F. Schulz for detailed exposition. To state the result in \cite{Sch},  we consider the Monge-Amp\`{e}re equation for a function $z=z(x,y)$ on a domain ${\mathcal{D}}\subset \mathbb{R}^2:$
\begin{equation}\label{50}
(z_{xx}+C)(z_{yy}+A)-(z_{xy}-B)^2=K(x,y,z)D(x,y,z,z_x,z_y)>0.
\end{equation}
where $A,B,C,D$ are functions of $x,y,z,p,q,$ and $p=z_x,q=z_y.$

\textbf{Assumption i)} $z\in C^{1,1}(\mathcal{D})$ and \begin{equation}\label{3.2}|z_x|+|z_y|\leq \mathcal{K}_1 .\end{equation}

\textbf{Assumption ii)} $A,B,C\in C^{1}(\mathcal{D}\times \mathbb{R}^3),
$ $K\in C^{\mu}(\mathcal{D}\times \mathbb{R}),$ for some $0<\mu<1,$
$D\in C^{1}(\mathcal{D}\times \mathbf{R}^3)$ and

\begin{equation}\label{3.3}
|A|+|B|+|C|+|D| \leq \mathcal{A}_1,\end{equation}
\begin{equation}\label{3.4} K,D\geq \frac{1}{\mathcal{A}_2},
\end{equation}
\begin{equation}\label{3.5}
 |\partial_{\mathcal{D}\times \mathbb{R}^3} A|+\cdots+
|\partial_{\mathcal{D}\times \mathbb{R}^3} D|\leq \mathcal{A}_3,\end{equation}
\begin{equation}\label{3.6} |K|_{C^{\mu}(\mathcal{D}\times \mathbb{R})}\leq \mathcal{A}_4 .
\end{equation}

 \textbf{Assumption iii)} The functions \begin{equation}\label{3.7}
\begin{split}
\phi_1(x,y)&=A_p,\\
\phi_2(x,y)&=A_q+2B_p,\\
\phi_3(x,y)&=C_p+2B_q,\\
\phi_4(x,y)&=C_q\\
\end{split}
\end{equation} are Lipschitz continuous with
\begin{equation}\label{3.8}
[\phi_1]_{0,1}^{\mathcal{D}}+\cdots+[\phi_4]_{0,1}^{\mathcal{D}}\leq \mathcal{A}_5.
\end{equation}

\begin{thm} \label{th3.1}(Theorem 9.4.1 in \cite{Sch})
Suppose $z\in C^{1,1}(\mathcal{D})$ is a solution of (\ref{50})
 such that the above Assumptions i), ii), iii) hold with the constants $\mathcal{K}_1,\mathcal{A}_1,\cdots, \mathcal{A}_5.$ Then $z\in C^{2,\mu}_{loc}(\mathcal{D})$, and for any $\mathcal{D}'\subset\subset \mathcal{D}$ there is an interior estimate
 \begin{equation}\label{3.9}
 \|\partial^2 z\|_{C^{\mu}(\mathcal{D}')}\leq C(\mu,\mathcal{K}_1,\mathcal{A}_1\cdots \mathcal{A}_5,dist(\mathcal{D}',\partial \mathcal{D})).
 \end{equation}
\end{thm}

For any  $\tilde{x}\in M,$ to invoke the result in \cite{Sch}, we fix a local coordinate system $(x,y)\in \mathcal{D}$ in $M$ around $\tilde{x}.$
Take $z(x,y)$ to be solution $u(x,y)$ of equation (\ref{Diri}) defined on $\Omega\supset \mathcal{D}.$
Then we find \begin{equation}\label{3.10}
\begin{split}
A&=g_{22}-\Gamma_{22}^{k}p_k,\\
B&=-g_{12}+\Gamma_{12}^{k}p_k,\\
C&=g_{11}+\Gamma_{11}^{k}p_k,\\
D&=g^{kl}p_kp_l+2z,\\
K(x,y,z)&=-K_g(x,y)det(g_{ij}),
\end{split}
\end{equation}
where $p_1=p,p_2=q.$

Note that by Propositions \ref{0order}, \ref{1order}, \ref{0'order}, we have estimated the upper bound of $u,$ $\nabla u$ and the lower bound of $u$ in the coordinate system $\mathcal{D},$ this gives rise to a control of the constants $\mathcal{K}_1,\mathcal{A}_1\cdots \mathcal{A}_5$
in terms of the geometry of $(\mathcal{D},g).$ From Theorem \ref{th3.1},  we have immediately
 \begin{prop} \label{prop3.2} For any nonnegative integer
$k\geq0,$ there exists a constant $D_{k}>0$ such that
\begin{equation}\label{3.11}\sup_{\Omega_{l}\supset B(x_0,k+1)}
|u_{l}|_{C^{2,\mu}(\bar{B}(x_0,k+\frac{1}{2}))} \leq D_{k},
\end{equation}
where the norm $C^{2,\mu}(\bar{B}(x_0,k+\frac{1}{2}))$ can be defined on  some (and
any) fixed finite coordinate covering of $\bar{B}(x_0,k+\frac{1}{2}).$
\end{prop}

We proceed to consider the third and higher order estimates
(\ref{est}). This may be done by standard Schauder estimate for elliptic equations.
\begin{prop}
\label{prop3.3} For any nonnegative integer
$k\geq0,$ there exists a constant $D_{k}>0$ such that
\begin{equation}\label{3.11}\sup_{\Omega_{l}\supset B(x_0,k+1)}
|u_{l}|_{C^{k}(\bar{B}(x_0,k))} \leq D_{k},
\end{equation}
where the norm $C^{k}(\bar{B}(x_0,k))$ can be defined on  some (and
any) fixed finite coordinate covering of $\bar{B}(x_0,k).$
\end{prop}
\begin{pf}
  By (\ref{Diri}), we see that $\nabla_iu$ satisfies an equation of the
following type
\begin{equation}\label{40}
\hat{g}^{jm}v_{jm}=f(x,v,\nabla v),
\end{equation}
where $\hat{g}=\nabla^2 u+g$.  By
previous second order estimate, we know (\ref{40}) is uniformly
elliptic on
${B}(x_0,k+\frac{1}{2})$
and the $C^{\mu}$ norm of $\hat{g}$ and $f$ are uniformly bounded (independent of $l$). The result follows from the standard interior Schauder estimate and bootstrap argument.
\end{pf}

\noindent {\sl \underline{Proof of Theorem \ref{t2.2}}}. By Proposition \ref{prop3.3} and Arzela-Ascoli theorem, we may extract a $C^{\infty}_{loc}$ convergent subsequence of $u_l.$ The limit is the desired solution. \hfill $\Box$

 Theorem \ref{t1} follows from Theorem \ref{t2.2} (see Remark \ref{rem1}).
\vskip0.3cm

\section{Proof of Theorems \ref{t1.2} and \ref{t1.3}}

The proof of Theorem \ref{t1.2} is based on Proposition
\ref{0order} and Theorem \ref{th3.1} and the following result.
\begin{prop}\label{p4.1}
   Under the assumptions (\ref{1.3}) (\ref{1.4}) of Theorem \ref{t1.2}, there exists $R>0$,
   such that $M$ admits a covering of isothermal coordinate charts $\{(U_i,(u^1,u^2))\}$
   with $U_i=\{ (u^1)^2+(u^2)^2<R^2\}$ such that\\
   1) for any $y_0\in M,$ there is an $
   U_{i_0}$ with $y_0\in \{(u^1)^2+(u^2)^2< \frac{R^2}{4}\}\subset U_{i_0}$,\\
   2) in each $U_i$, the metric $g$ of $M$ takes the form $$g=\psi[(du^1)^2+(du^2)^2],$$
   and it satisfies  \begin{equation}\label{43}
         \begin{split}
   &c^{-1}\leq \psi\leq c,\\
    &  |\psi|_{C^{2,\mu}(U_{i})}\leq c_{\mu}, \ \ \ \ \
      \end{split}
   \end{equation}
   where $c$ and $c_{\mu}$ are constants independent of
   $i$. Moreover, if the additional (\ref{42}) is satisfied, we have
\begin{equation}\label{44}
         \begin{split}
      |\psi|_{C^{l+1,\alpha}(U_{i})}\leq c_{l,\alpha} \ \ \ \ \ \text{for any } \alpha\in
    (0,1),
          \end{split}
   \end{equation}
where $c_{l,\alpha}$ are constants independent of $i.$
\end{prop}

 \noindent \underline{\sl Proof of Theorem \ref{t1.2}}.  By Proposition
\ref{0order} and Theorem \ref{t1}, we know there exists a smooth
isometric embedding $X:M\rightarrow
 \mathbb{R}^{2,1} $ such that  $u=-\frac{1}{2}\langle X,X\rangle$ satisfying
$\frac{1}{2{C}_2}\leq u\leq \frac{1}{2C_1}. $  Let $R$ be the constant provided in Proposition \ref{p4.1}. Let the coordinates $(x,y)$ in equation (\ref{50}) to be the isothermal coordinates $(u^1,u^2)$ in Proposition \ref{p4.1}, $z(x,y)=u(x,y),$ and $\mathcal{D}=\{x^2+y^2<\frac{R^2}{4}\}$. In these coordinates, (\ref{3.10}) becomes  \begin{equation}
\begin{split}
A&=\psi-\Gamma_{22}^{k}p_k,\\
B&=\Gamma_{12}^{k}p_k,\\
C&=\psi+\Gamma_{11}^{k}p_k,\\
D&=\psi^{-1}(p_1^2+p_2^2)+2z,\\
K(x,y,z)&=-K_g(x,y)\psi^2.
\end{split}
\end{equation}

Estimate (\ref{43}) and the Proposition \ref{1order} imply that there is a constant $C$ depending only on $C_1,C_2,C_{\mu}$ such that the constants in (\ref{3.2})-(\ref{3.8}) can be bounded by $C$
$$
\mathcal{K}_1,\mathcal{A}_1,\cdots,\mathcal{A}_5\leq C.
$$
Theorem \ref{th3.1} implies $|\partial_{ij}u|_{C^{\mu}(B(0,\frac{R}{4}))}\leq C.$
Combining this with (\ref{43}) particularly gives $|h_{ij}|\leq C.$ This proves ii) in Theorem \ref{t1.2}.
\hfill$\Box$
\vskip 0.3cm
If the additional higher covariant derivatives bound (\ref{42}) is assumed, notice that (\ref{44}) holds, then (\ref{010}) follows by the same argument in  Proposition \ref{prop3.3}.

\vskip 0.4cm
 \noindent \underline{\sl Proof of Theorem \ref{t1.3}}. After an isometry $\tilde{\iota}$ of $\mathbb{R}^{2,1},$ $\tilde{\iota}\circ \tilde{X}(M)$ can be pinched between the light cone and a hyperboloid associated to $X,$ and we can define
$\tilde{u}=-\frac{1}{2}\langle  \tilde{\iota}\circ \tilde{X},\tilde{\iota}\circ \tilde{X}\rangle,$ which satisfies
 $ 0 < \tilde{u}\leq C. $

 Using the polar coordinates in Remark \ref{rem1},
we know that $\tilde{\iota}\circ\tilde{X}$ is determined by $\tilde{u}$ and an
isometry
$\tilde{i}:(M,\frac{g+(d\sqrt{2\tilde{u}})^2}{2\tilde{u}})\rightarrow
\mathbb{H}.$  To show that $\tilde{\iota}\circ\tilde{X}$ is congruent to  $X,$ it
suffices to show that $u=\tilde{u}.$ Indeed, once we have
$u=\tilde{u},$ then $\tilde{\iota}\circ\tilde{X}=\sigma \circ X,$ where
$\sigma=\tilde{i}i^{-1}\in Iso(\mathbb{H})\subset
Iso(\mathbb{R}^{2,1}).$ Then $\tilde{X}=\iota\circ X,$ where $\iota=\tilde{\iota}^{-1}\circ \sigma.$

We need some a prior estimates of $\tilde{u}$ up to second order.
 The powerful tool is  the  Cheng-Yau's maximum
principle \cite{CY}, since the curvature is assumed to be bounded:   for any $C^2$ function $F$
bounded from above, there is a sequence of $x_k\in M$ and
$\varepsilon_{k}\rightarrow 0$ such that
\begin{equation}\label{4.6}
\begin{split}
a)\ &\sup_{M}F-F(x_k)\leq \varepsilon_k,\\
b)\ & |\nabla F|(x_k)\leq \varepsilon_k,\\
c)\  &\nabla^2 F(x_k)\leq \varepsilon_k g.
\end{split}
\end{equation}

Applying the above maximum principle to $\tilde{u}$ and $-\tilde{u},$ we immediately know $$\frac{1}{2C_2}\leq \tilde{u}\leq \frac{1}{2C_1}.$$

 We claim that the gradient of $\tilde{u}$ is also bounded, and more precisely, it satisfies $$|\nabla \tilde{u}|\leq \frac{1}{\sqrt{C_1}}.$$

 Indeed, for any
$\tilde{x}\in M,$ let $\gamma$ be a geodesic of unit speed such that $\gamma(0)=\tilde{x}.$
we would like to control $|\frac{d}{dt}(\tilde{u}\circ \gamma)(0)|.$ By the  convexity of the function
$\tilde{u}+\frac{1}{2}d^2(\tilde{x},\cdot),$ we know
$$
t |\frac{d}{dt}(\tilde{u}\circ \gamma)(0)|\leq \max\{\tilde{u}(\gamma(t))-\tilde{u}(\gamma(0)),\tilde{u}(\gamma(-t))-\tilde{u}(\gamma(0)) \}+\frac{t^2}{2}
$$
It follows that $|\frac{d}{dt}(\tilde{u}\circ \gamma)(0)|\leq \frac{1}{\sqrt{C_1}}$ by taking $t=\frac{1}{\sqrt{C_1}}.$ This implies
$|\nabla \tilde{u}|(\tilde{x})\leq \frac{1}{\sqrt{C_1}},$ and the claim is proved.

Combining the gradient estimate of $\tilde{u}$ with the proof of Proposition \ref{prop3.2}, we know
that $|\nabla^2\tilde{u}|$ is bounded.

Summarizing the above estimates, it follows that there is $C>0$ such that
\begin{equation}\label{uni}
 \frac{1}{C}\leq u_t\leq C, \ \ |\nabla u_t|\leq C, \ \nabla^2{u}_t+g\geq C^{-1}g,
\end{equation}
where $u_{t}=u+t(\tilde{u}-u),$ $t\in [0,1].$





Note that $u$ and $\tilde{u}$ satisfy the same equation (\ref{4}),
this implies
\begin{equation}\label{EQ}
\int_{0}^{1}(g+\nabla^2u_t)^{ij}dt\nabla_{ij}^2(\tilde{u}-u)=\int_{0}^{1}\frac{2\langle
\nabla u_t,\nabla(\tilde{u}-u)\rangle+2(\tilde{u}-u)}{|\nabla
u_t|^2+2u_t}dt
\end{equation}
Let $F=\tilde{u}-u$ in (\ref{4.6}), combining (\ref{4.6}) with (\ref{uni})
(\ref{EQ}), we have
$$C\varepsilon_{k}\geq \sup_{M}(\tilde{u}-u).
$$
This gives $u\geq \tilde{u}.$ Similarly, we have $u\leq \tilde{u}.$
Hence $u=\tilde{u}.$

To prove ii), one can show $u\circ \gamma=u$  for any
$\gamma\in Iso(M,g)$ by the above Cheng-Yau's maximum principle (\ref{4.6}).
This implies $Iso(M,g)\subset Iso(M,\bar{g}),$
where $\bar{g}=\frac{g+(d\sqrt{2u})^2}{2u}.$ The injective
homomorphism $\rho: Iso(M,g)\rightarrow Iso(\mathbb{H})$ is given by
$\rho(\gamma)=i \circ \gamma\circ i^{-1}.$
\hfill$\Box$
\vskip0.3cm

\section{Appendix}
The purpose of this appendix is to give an alternative method for
second order estimate. The argument we present here is classical,
straightforward,  and may be generalized to higher dimensions (see
\cite{GL1}). The price to be paid is that this method requires some
geometry of the background manifold. It works well on those points
where the values of $u$ (solution of (\ref{Diri})) are not too
large, comparing to the local geometry.
\begin{prop}\label{prop5.1} There exists ${C}>0$
depending only on $C_1$ satisfying the following property. Fix
$\tilde{x}\in M,$ suppose there exist a real number $r_0>0$ and  a
solution $u$ of (\ref{Diri}) defined on domain  $\Omega_{l}\supset
B(\tilde{x},r_0)$ satisfying

\begin{equation}\label{assump}
u(\tilde{x})< \frac{r_0}{2\sqrt{\max\limits_{y\in
\bar{B}(\tilde{x},r_0)}(-K_g(y))}\coth (\sqrt{\max\limits_{y\in
\bar{B}(\tilde{x},r_0)}(-K_g(y))}r_0)}.
\end{equation}
Then
\begin{equation}
\begin{split}
(g+\nabla^2u)({x}) \leq & \frac{e^{Cc_2'}}{
\frac{r_0}{2\sqrt{c_2}\coth (\sqrt{c_2}r_0)}
-u(\tilde{x})}(1+\sqrt{c_4}\frac{r_0}{\sqrt{c_2}}+c_2'(1+\frac{r_0}{\sqrt{c_2}}+c_3\frac{r_0}{\sqrt{c_2}})),
\end{split}
\end{equation}
on
${B(\tilde{x},\frac{\sqrt{C_1}(\frac{r_0}{2\sqrt{c_2}\coth(\sqrt{c_2}
r_0)}-u(\tilde{x}))}{6})},$ where
\begin{equation}
\begin{split}
 c_2&=\max_{x\in \bar{B}(\tilde{x},r_0)}(-K_g(x)),\ \ \ \ \ \ \ \ \ \
 \
 c_2'=\max_{x\in \bar{B}(\tilde{x},r_0+1)}(-K_g(x)),\\
 c_3&=\max_{x\in \bar{B}(\tilde{x},r_0)}| \nabla\log
(-K_g(x))|, \ \ c_4=\max_{x\in \bar{B}(\tilde{x},r_0)}|
\nabla^2\log(- K_g(x))|.
\end{split}
\end{equation}
\end{prop}
\vskip 0.3cm

Note that by Proposition \ref{0order},  condition (\ref{assump}) can be justified at
each  $\tilde{x}$ (for suitable $r_0$) when the curvature $K$ satisfies
 \begin{equation}\label{-1}-{C}_2^2 (d(x,x_0)+{C}_3)^2 \leq K(x)\leq -C_1
 \end{equation}
 for some $x_0\in M$ and positive constants $0<{C}_2<C_1<{C}_3.$

\noindent  {\sl \underline{Proof of Proposition \ref{prop5.1}}.}
Consider an auxiliary function:
\begin{equation}\label{q}\left\{
\begin{split} &STM\rightarrow \mathbb{R}\\&(x,\gamma)\rightarrow \eta(x)
(1+\nabla_{\gamma\gamma}u)e^{\frac{a}{2}(|\nabla u|^2+2u)(x)},
\end{split}
\right. \end{equation} where $x\in M,\gamma\in T_xM,|\gamma|=1,$
$STM$ is the unit tangent bundle  of $M$, $\eta$ is a cutoff
function on $M$ and $a\geq 1$ is a constant to be specified  later.
Suppose the maximum
\begin{equation*} \max_{(x,\gamma)\in STM}\eta(1+\nabla_{\gamma\gamma}u)e^{\frac{a}{2}(|\nabla
u|^2+2u)} \end{equation*}
 is achieved at $\bar{x}\in $ $\text{supp}(\eta)$ for
some $\gamma\in T_{\bar{x}}M$ with $|\gamma|=1.$ Diagonalize
$u_{ij}=\lambda_i \delta_{ij}$ at $\bar{x}$ with the orthonormal
eigenvectors $e_i.$ Let $e_1=\gamma.$ Parallel transport each $e_i$
along radial geodesics, we obtain a field of  orthonormal  frame
$\{e_i\}$ near $\bar{x}.$
The function $w=\eta(1+\nabla_{e_1,e_1}u)e^{\frac{a}{2}(|\nabla
u|^2+2u)}$ defined near $\bar{x}$ achieves
its maximum at $\bar{x}.$ In the following, we use $C$ to denote various big constants
depending only on $C_1.$

 At the point $\bar{x},$ we have
\begin{equation}\label{116}
   0=\nabla_i\log
  w=\frac{\nabla_{i11}u}{1+\lambda_1}+a(1+\lambda_i)u_i+\frac{\nabla_i\eta}{\eta},
\end{equation}
\begin{equation}\label{17}
\begin{split}
 0\geq\nabla_{ij}\log w=&\frac{\nabla_{ij11}u}{1+\lambda_1}-\frac{\nabla_{i11}u\nabla_{j11}u}{(1+\lambda_1)^2}
 +a(u_k\nabla_{ijk}u+(\lambda_i+\lambda_i^2)\delta_{ij})\\&
+\frac{\nabla_{ij}\eta}{\eta}-\frac{\nabla_i\eta\nabla_j\eta}{\eta^2}.
\end{split}
\end{equation}

Let $f=f(x,z,p)\triangleq \log(-K)+\log(|\nabla u|^2+2u)),$ where $z=u,p=\nabla u.$  Differentiating
the equation (\ref{Diri}), we get
\begin{equation}\label{23}
\frac{\nabla_{kii}u}{1+\lambda_i}=\nabla_k f,
\end{equation}
\begin{equation}\label{24}
\frac{\nabla_{11ii}u}{1+\lambda_i}-\frac{(\nabla_{1ij}u)^2}{(1+\lambda_i)(1+\lambda_j)}=f_{11}.
\end{equation}
Combining (\ref{17}), (\ref{23}), (\ref{24}) and Ricci formula, we have
\begin{equation}\label{32}
\begin{split}
&(1+\lambda_1)(-\frac{\nabla_{ii}\eta}{(1+\lambda_i)\eta}+\frac{|\nabla_i\eta|^2}{(1+\lambda_i)\eta^2})\\
&+\frac{1}{1+\lambda_i}[\frac{(\nabla_{i11}
 u)^2}{1+\lambda_1}-a(R_{ijip}u_pu_j+\lambda_i+\lambda_i^2)(1+\lambda_1)\\
 &-(\nabla_iR_{i11p}+\nabla_1R_{i1ip})u_p-2R_{1i1i}(\lambda_1-\lambda_i)]
  -\frac{(\nabla_{1ij}u)^2}{(1+\lambda_i)(1+\lambda_j)}\\ \geq & \  f_{11}+a
u_j\nabla_jf(1+\lambda_1).
  \end{split}
\end{equation}
By direct computations, we have (at $\bar{x}$)

\begin{equation}\label{33}
\begin{split}
f_{11}+a u_k\nabla_kf(1+\lambda_1)\geq &\frac{2}{|\nabla
u|^2+2u}\langle \nabla u, -a\nabla u-\frac{\nabla\eta}{\eta}\rangle (1+\lambda_1)\\
&+\log(-K)_{11}+a\langle \nabla\log(-K),\nabla u\rangle(1+\lambda_1)\\
&+\frac{a(1+\lambda_1)|\nabla u|^2}{|\nabla
u|^2+2u}-\frac{8(|\lambda_1|+1)^2}{|\nabla
u|^2+2u}+2\frac{R_{1j1l}u_ju_l}{|\nabla u|^2+2u}.
  \end{split}
\end{equation}
By (\ref{116}), we have
\begin{equation}\label{34}
\begin{split}
\sum_{i}\frac{|\nabla_i\eta|^2}{(1+\lambda_i)\eta^2}(1+\lambda_1)
&=\frac{|\nabla_1\eta|^2}{\eta^2}+\sum_{i\geq
2}\frac{|\nabla_{i11}u|^2}{(1+\lambda_1) (1+\lambda_i)}\\& \ \ \
-2a(1+\lambda_1)\nabla_iu\frac{\nabla_i\eta}{\eta}
-a^2(1+\lambda_1)(1+\lambda_i)u_i^2.
\end{split}
\end{equation}
Note that
\begin{equation}\label{35}
\begin{split}
&\sum_{i,j}\frac{(\nabla_{1ij}u)^2}{(1+\lambda_i)(1+\lambda_j)}-\sum_{i}
\frac{(\nabla_{i11}u)^2}{(1+\lambda_1)(1+\lambda_i)}-\sum_{i\geq
2}\frac{|\nabla_{i11}u|^2}{(1+\lambda_1) (1+\lambda_i)}\\
\geq & -4\frac{|\nabla u|}{|\nabla u|^2+2u}\frac{|\nabla
\eta|}{\eta}(1+\lambda_1),
\end{split}
\end{equation}
\begin{equation}\label{36}
\begin{split}
&\sum_{i}-\frac{1}{1+\lambda_i}(R_{ijip}u_pu_j+\lambda_i+\lambda_i^2)(1+\lambda_1)\\
\leq & -\frac{2u}{|\nabla u|^2+2u}(1+\lambda_1)^2+2(1+\lambda_1),
\end{split}
\end{equation}
and
\begin{equation}\label{37}
\begin{split}
&\frac{1}{1+\lambda_i}[
  -(\nabla_iR_{i11p}+\nabla_1R_{i1ip})u_p-2R_{1i1i}(\lambda_1-\lambda_i)]\\
  \leq &
  \frac{2|\nabla u|}{|\nabla u|^2+2u}|\nabla\log(- K)|(1+\lambda_1)+\frac{2}{|\nabla
  u|^2+2u}(1+\lambda_1)^2+2K.
  \end{split}
\end{equation}
Multiplying both sides of (\ref{32}) by $\eta^2$, combining
(\ref{33})-(\ref{37}), we get

\begin{equation}\label{38}
\begin{split}
&L_1(1+\lambda_1)^2\eta^2-L_2(1+\lambda_1)\eta-L_3\leq
\eta(1+\lambda_1)\sum_{i\geq 1}\frac{-\nabla_{ii}\eta}{1+\lambda_i},
  \end{split}
\end{equation}
where
\begin{equation}\label{L}
\begin{split}
L_1&=a\frac{2u}{|\nabla u|^2+2u}-\frac{10}{|\nabla u|^2+2u},\\
L_2&= (6\frac{|\nabla u|}{|\nabla u|^2+2u}+2a|\nabla u|)|\nabla
\eta|+2a\eta+\frac{|\nabla u|^2}{|\nabla u|^2+2u}a\eta\\&\ \ \ \
+a|\nabla \log(-K)||\nabla u|\eta+\frac{2|\nabla u|}{|\nabla
u|^2+2u}|\nabla\log(-K)|\eta,\\
L_3&=|\nabla \eta|^2+\eta^2|\nabla^2\log(-K)|+2K \frac{|\nabla_1
u|^2+2u}{|\nabla u|^2+2u}\eta^2.
\end{split}
\end{equation}

 Note that by (\ref{assump}), Lemma \ref{39} is applicable. Choose the cutoff function
$\eta$ in (\ref{q}) to be  $\varphi^{\tilde{x}}$ in Lemma \ref{39}, and
consider the maximum of the quantity $w$ on $Q_{\tilde{x}}.$ From
iii) in Lemma \ref{39}, we have
\begin{equation}\label{61}
\eta(1+\lambda_1)\sum_{i\geq
1}\frac{-\nabla_{ii}\eta}{1+\lambda_i}\leq 2(1+\lambda_1)\eta.
\end{equation}

 Since $u(\bar{x})\geq C^{-1}{c'_2}^{-1}$ by
Corollary \ref{0''order},  choosing $a={Cc_2'}$ in (\ref{L}) and
combining Lemma \ref{39}, we have
\begin{equation}\label{60}
\begin{split}
L_1&\geq \frac{10}{|\nabla u|^2+2u}\geq 2C_1,\\
L_2&\leq c_2'(1+\frac{r_0}{\sqrt{c_2}}+c_3\frac{r_0}{\sqrt{c_2}}),\\
L_3&\leq C(1+c_4)\frac{r^2_0}{{c_2}}.
\end{split}
\end{equation}
 From (\ref{38}), (\ref{60}), (\ref{61}), we have
\begin{equation}\label{70}
\begin{split}
(1+\lambda_1)\eta &\leq
\max\{\sqrt{\frac{2L_3}{L_1}},\frac{2(L_2+2)}{L_1}\}\\
&\leq
C(1+\sqrt{c_4}\frac{r_0}{\sqrt{c_2}}+c_2'(1+\frac{r_0}{\sqrt{c_2}}+c_3\frac{r_0}{\sqrt{c_2}})).
\end{split}
\end{equation}

Combining Lemma \ref{39} i) and (\ref{70}),
we have \begin{equation}
\begin{split} (1+\lambda_1)(x)\leq &
\frac{e^{Cc_2'}}{ \frac{r_0}{2\sqrt{c_2}\coth(\sqrt{c_2}r_0)}
-u(\tilde{x})}\\ & \times (1+\sqrt{c_4}\frac{r_0}{\sqrt{c_2}}
+c_2'(1+\frac{r_0}{\sqrt{c_2}}+c_3\frac{r_0}{\sqrt{c_2}}))
\end{split}
\end{equation}
on
${B(\tilde{x},\frac{\sqrt{C_1}(\frac{r_0}{2\sqrt{c_2}\coth(\sqrt{c_2}r_0)}-u(\tilde{x}))}{6})}.$



The proof of Proposition \ref{prop5.1} is completed.\hfill $\Box$

\begin{rem} The most computations in this section are just modifications of those in the classical
theory of Monge-Amp\`{e}re equations. A closer reference is
\cite{GL1}, where the Dirichlet problem of real
Monge-Amp\`{e}re equations on manifolds has been
systematically studied. The observation is that these
estimates could be localized under certain geometric conditions.
\end{rem}

Bing-Long Chen \\
 Department of Mathematics, \\
Sun Yat-Sen University,\\
Guangzhou, P.R.China, 510275\\
Email: mcscbl@mail.sysu.edu.cn\\

\noindent Le Yin\\
College of Mathematics and Computational Science,\\
Shenzhen University,\\
   Shenzhen, P.R.China, 518060\\
    Email: lyin@szu.edu.cn
\end{document}